\documentclass[10pt,draft]{amsart}
\usepackage{amsmath,amssymb,amsthm}

\begin{document}
\newtheorem{lem}{Lemma}[section]
\newtheorem{prop}{Proposition}[section]
\newtheorem{cor}{Corollary}[section]
\numberwithin{equation}{section}
\newtheorem{thm}{Theorem}[section]
\theoremstyle{remark}
\newtheorem{example}{Example}[section]
\newtheorem*{ack}{Acknowledgment}
\theoremstyle{definition}
\newtheorem{definition}{Definition}[section]
\theoremstyle{remark}
\newtheorem*{notation}{Notation}
\theoremstyle{remark}
\newtheorem{remark}{Remark}[section]
\newenvironment{Abstract}
{\begin{center}\textbf{\footnotesize{Abstract}}%
\end{center} \begin{quote}\begin{footnotesize}}
{\end{footnotesize}\end{quote}\bigskip}
\newenvironment{nome}
{\begin{center}\textbf{{}}%
\end{center} \begin{quote}\end{quote}\bigskip}

\newcommand{\triple}[1]{{|\!|\!|#1|\!|\!|}}
\newcommand{\xx}{\langle x\rangle}
\newcommand{\ep}{\varepsilon}
\newcommand{\al}{\alpha}
\newcommand{\be}{\beta}
\newcommand{\de}{\partial}
\newcommand{\la}{\lambda}
\newcommand{\La}{\Lambda}
\newcommand{\ga}{\gamma}
\newcommand{\del}{\delta}
\newcommand{\Del}{\Delta}
\newcommand{\sig}{\sigma}
\newcommand{\ome}{\omega}
\newcommand{\Ome}{\Omega}
\newcommand{\C}{{\mathbb C}}
\newcommand{\N}{{\mathbb N}}
\newcommand{\Z}{{\mathbb Z}}
\newcommand{\R}{{\mathbb R}}
\newcommand{\Rn}{{\mathbb R}^{n}}
\newcommand{\Rnu}{{\mathbb R}^{n+1}_{+}}
\newcommand{\Cn}{{\mathbb C}^{n}}
\newcommand{\spt}{\,\mathrm{supp}\,}
\newcommand{\Lin}{\mathcal{L}}
\newcommand{\SSS}{\mathcal{S}}
\newcommand{\F}{\mathcal{F}}
\newcommand{\eei}{\langle\eta\rangle}
\newcommand{\xei}{\langle\xi-\eta\rangle}
\newcommand{\yy}{\langle y\rangle}
\newcommand{\dint}{\int\!\!\int}
\newcommand{\hatp}{\widehat\psi}
\renewcommand{\Re}{\;\mathrm{Re}\;}
\renewcommand{\Im}{\;\mathrm{Im}\;}

\title[ Counterexamples of Strichartz inequalities ]%
{Counterexamples of Strichartz inequalities for Schr\"odinger equations with 
repulsive potentials
}

\author{Michael Goldberg}

\address{Michael Goldberg\\Department of Mathematics, Johns Hopkins University\\
Baltimore, MD 21218, USA}

\email{mikeg@math.jhu.edu}

\author{Luis Vega}

\address{Luis Vega\\
Universidad del Pais Vasco, Apdo. 64\\
48080 Bilbao, Spain}

\email{mtpvegol@lg.ehu.es}

\author{Nicola Visciglia}

\address{Nicola Visciglia\\
Dipartimento di Matematica Universit\`a di Pisa\\
Largo B. Pontecorvo 5, 56100 Pisa, Italy}

\email{viscigli@mail.dm.unipi.it}




\thanks{\noindent This research was supported by HYKE (HPRN-CT-2002-00282).
The second author was supported also by a MAC grant (MTM 2004-03029) 
and the third one by an INDAM 
(Istituto Nazionale di Alta Matematica) fellowship}

\date{}

\begin{abstract} In each dimension $n \ge 2$, we construct a class of 
nonnegative potentials that are homogeneous of order $-\sigma$, chosen from
the range $0 \le \sigma < 2$, and 
for which the perturbed Schr\"odinger equation does not satisfy
global in time Strichartz estimates.

\end{abstract}

\maketitle
\vspace{0.2cm}

\section{Introduction}

The study of dispersive estimates 
for evolution equations has received a lot of attention in the literature.
\vspace{0.1cm}

In particular the free Schr\"odinger equation
\begin{equation}\label{eq.freesch}
\begin{cases}
{\bf i}\partial_{t}u-\Delta_x u=0, 
(t, x)\in {\mathbf R}_t \times {\mathbf R}_x^n,\\
 u(0,x)=f(x),
\end{cases}
\end{equation}
exhibits a rich set of dispersive and smoothing estimates. 

\vspace{0.1cm}

An important class of estimates satisfied by the solutions of
the Cauchy problem \eqref{eq.freesch} are the Strichartz estimates
that we recall below:
\begin{equation}\label{strichartz}
\|u(t, x)\|_{L^p({\mathbf R}; L^q({\mathbf R}^n))} 
\leq C \|f\|_{L^2({\mathbf R}^n)}
\end{equation}
provided that
\begin{equation}\label{gap}
\frac 2p + \frac nq=\frac n2, \hspace{0.2cm} 2\leq p\leq \infty, 
\hspace{0.2cm} (n, p)\neq (2, 2),
\end{equation}
where $u(t, x)$ is the unique solution of \eqref{eq.freesch},
see \cite{KT}.

\vspace{0.1cm}

Notice that estimates \eqref{strichartz} 
describe a certain regularity for the solutions of \eqref{eq.freesch}
in terms of summability but they do not give 
any gain of derivatives.
For this reason Strichartz estimates 
are very useful in order to treat the local and global well-posedness
of the semilinear Schr\"odinger equation,
but they are useless
in the study of the local and global well-posedness 
of the nonlinear Schr\"odinger equation
with nonlinearities which involve derivatives.

\vspace{0.1cm}

In this case the following local smoothing 
estimate has turned out to be very useful
(see \cite{CS}, \cite{Sj}, \cite{V}):
\begin{equation}\label{smoothing}
\sup_{R\in (0, \infty)} \frac 1R \int_{-\infty}^\infty 
\int_{B_R} |\nabla_x u|^2 dxdt
\leq C \|f\|_{\dot H^\frac 12({\mathbf R}^n)},
\end{equation}
where $\dot H^s({\mathbf R}^n)$ denote 
the usual homogeneous Sobolev spaces.
In fact
the reverse inequality is almost true,
see \cite{VV}.

\vspace{0.2cm}

It is a natural question to understand if the estimates \eqref{strichartz}
and \eqref{smoothing} extend when a potential
 is added to \eqref{eq.freesch}, i.e. 
\begin{equation}\label{CP}
  \begin{cases}  {\bf i} \partial_{t}u-\Delta_x u+V(t, x)u=0,
 (t, x)\in {\mathbf R}_t \times {\mathbf R}^n_x,\\
 u(0, x)=f(x).
\end{cases}
\end{equation}

There is a huge literature on the subject, let us recall 
in particular \cite{BRV}, \cite{BPST1}, \cite{BPST2}, 
\cite{CoSa}, \cite{DPV}, \cite{G}, \cite{JSS}, \cite{RS}, 
\cite{RV1}, \cite{RV2}, \cite{S}, \cite{ST},
even if the list is very far from being complete.
It is important to observe that 
the conditions on the potential change dramatically  depending on 
whether or not the estimates are global or local in time. 
Our interest in this paper will be on global in time inequalities. 
At this respect in \cite{BRV} the authors have been able 
to prove that
the local smoothing estimate \eqref{smoothing} 
is still satisfied provided that
the free Schr\"odinger equation is 
perturbed by a repulsive potential. By repulsive it is meant
that the positive part of 
the radial derivative is small in an appropriate sense.

\vspace{0.1cm}

In this paper we shall focus on the Cauchy problem \eqref{CP}
in the case that $V(x)$ is a potential which is homogeneous of order $-\sigma$,
i.e. 
\begin{equation}
\label{hom}V(\lambda x)=\lambda^{-\sigma}V(x), \forall (\lambda, x)
\in (0, \infty)\times {\mathbf R}_x^n.
\end{equation} 

Therefore
$$V(x)=|x|^{-\sigma}V\left( \frac x{|x|}\right), 
\forall x\in {\mathbf R}^n_x \setminus\{0\},$$
and 
$V(x)$ is uniquely determined 
by its restriction on the sphere
$${\mathbf S}^{n-1}=\{(x_1,..., x_n)\in {\mathbf R}^n_x
\hbox{ s.t. }  \sum_{i=1}^n x_i^2=1\}.
$$
So long as $V(x) \ge 0$ and $\sigma\ge 0$, such a potential will be repulsive,
as the radial derivatives are all non-positive.

\vspace{0.1cm}

The study of  this type of perturbations was started 
by H. Herbst in \cite{H} where the role
of the critical points of $V$ is emphasized. 
More recently in \cite{PV1} and \cite{PV2}
the authors proved for the reduced
wave (Helmholtz) equation an energy estimate that 
suggests concentration should occur 
in the set of critical points. Also in \cite{HS} 
the authors study the existence 
and completeness of the wave operator for  the same kind of potentials.

\vspace{0.1cm}

The aim of our paper is to analyze
whether the Strichartz estimates \eqref{strichartz} are
preserved by this class of perturbations. 
At this respect it is important to 
mention \cite{JSS} and \cite{ST} where the local smoothing 
estimate is used to conclude estimates as \eqref{strichartz}.
We will prove that given any $n \ge 2$ and $0\le \sigma < 2$ there is a class of
potentials satisfying \eqref{hom} for which \eqref{strichartz}
does not hold even though \eqref{smoothing} does.

\vspace{0.1cm}

In order to state our main result we need the following
definition.

\begin{definition}\label{def}
A function $V(x)$ that is homogeneous
of order $-\sigma$
is said to be of generalized Morse type provided that 
its restriction on the sphere
$$V_{|{\mathbf S}^{n-1}}:S^{n-1}\ni \omega
\rightarrow V(\omega)\in \mathbf R$$
has a nondegenerate minimum point, i.e. its Hessian at a minimum point
is a nondegenerate bilinear form. 
\end{definition}

\vspace{0.1cm}

We can now state the main result of this paper.

\vspace{0.1cm}

\begin{thm}\label{main}
Assume that $n\geq 2$ and 
$V\in C^3({\mathbf R}^n_x\setminus \{0\})$ 
is a homogeneous function of order $-\sigma$,\ $0 \le \sigma < 2$, 
which is of generalized Morse type.  Further assume that the minimum value 
of $V(x)$ on each sphere $\lambda {\mathbf S}^{n-1}$ is exactly zero.
Then for any $p\neq \infty$
the Strichartz estimates \eqref{strichartz} cannot be satisfied 
by the corresponding solutions of  \eqref{CP}.
\end{thm}

\begin{remark}
It is well known that Strichartz estimates are valid for the
potential $V(x) \equiv 0$, which satisfies every hypothesis of the
theorem except that it is not a generalized Morse function.
\end{remark}


\begin{remark}
Looking 
at the proof of
theorem \ref{main} it is easy to see that the
conclusion is still true 
in the case that 
\eqref{hom} is satisfied only for 
$\lambda, x$ large enough, 
assuming that 
the restriction of $V(x)$ on a sphere large enough is a 
generalized Morse function with minimum value zero.
\end{remark}

\begin{remark}
If the potential is homogeneous of degree zero the condition that the minimum has to be zero can be trivially relaxed by adding a real constant  to the potential, so that the corresponding solutions are changed just by a factor of modulus one that do not affect the validity of Strichartz estimates.
\end{remark}

Next we fix some notations useful in the sequel.

\vspace{0.2cm}

\begin{notation}
We shall denote by $x\in {\mathbf R}^n_x$ the full space variables.
In some cases we shall use the following splitting:
$$x=(y, z)\in {\mathbf R}^{n-1}_y \times {\mathbf R}_z.$$

Notice that this decomposition
allows us to split the full Laplacian operator $\Delta_x$ 
in the following way:
$$\Delta_x=\Delta_y + \partial_z^2.$$

For any $1\leq p, q<\infty$
and for any time-independent and time-dependent functions
$f(x)$ and $F(t,x)$, we 
shall use the following norms:
$$\|f(x)\|_{L^q_x}^q=\int_{{\mathbf R}^n} |f(x)|^q dx
\hspace{0.2cm} \hbox{ and } \hspace{0.2cm}
\|F(t, x)\|_{L^p_t L^q_x}^p=\int_{0}^\infty \|F(t, .)\|_{L^q_x}^p dt.$$

If $0<T<\infty$ then we shall write
$$\|F(t, x)\|_{L^p_T L^q_x}^p=\int_{0}^T \|F(t, .)\|_{L^q_x}^p dt.$$
and in some cases 
$$\|F(t, x)\|_{L^p_T L^q_x}= \|F(t, x)\|_{L^p((0,T); L^q_x)}.$$
\end{notation}

For any $1\leq p \leq \infty$ we shall denote by $p'$
the unique number such that
$$\frac 1p + \frac 1{p'}=1$$
and $1\leq p' \leq \infty$.

\vspace{0.4cm}

The work is organized as follows.
In section \ref{preliminary}
we make some preliminary computations
that will be useful in section \ref{sectionthm}
where we shall prove theorem \ref{main}.

\section{Preliminary Computations}\label{preliminary}

It is well-known
that for any symmetric positive definite matrix 
$\left ( a_{ij}\right )_{i,j=1,...,n-1}$
the operator 
\begin{equation}\label{schrop}-\Delta_y  + \sum_{i.j=1}^{n-1} a_{ij}y_i y_j,
\hbox{ }y\in {\mathbf R}_y^{n-1},  
\end{equation}
has compact resolvent
and in particular
its spectrum reduces to its point spectrum.
As a consequence of this fact
there exists a sequence of eigenvalues and eigenfunctions
associated to \eqref{schrop}.
In fact using a linear transformation
these operators reduce to the harmonic oscillator whose 
spectrum is explicitly 
well-known.

\vspace{0.1cm}

Therefore in the sequel we will assume that there 
exists a couple $(\lambda, v(y))$
such that 
\begin{equation}\label{eigen}-\Delta_y v +
\left ( \sum_{i,j=1}^{n-1} 
a_{ij} y_iy_j\right)v = \lambda v, \ \forall y\in \mathbf R^{n-1}_y
\end{equation}
and moreover $v(y)\in \cap_{p=1}^\infty L^p({\mathbf R}^{n-1}_y)$.

This choice will depend 
on the matrix $\left (a_{ij}\right )_{i,j=1,...,n-1}$
which in turn depends on the minimum of the potential.
Hereafter the couple $(\lambda, v(y))$ is fixed. 

\vspace{0.1cm}

Next we can introduce the rescaled function
\begin{equation}\label{w}w (y, z):=v \left ( \frac {y}{\sqrt{z^\beta}}\right), 
\hbox{  } \forall (y, z)\in {\mathbf R}^{n-1}_y\times (0, \infty),
\end{equation}
\noindent where $\beta = 1 + \frac{\sigma}2$ and $\sigma$ is the constant that appears
in theorem 
\ref{main}.  
By an elementary computation, this 
satisfies
\begin{equation}\label{rescaling}
-\Delta_y w + \frac{1}{z^{2\beta}}\left (\sum_{i,j=1}^{n-1} a_{ij}y_iy_j\right ) w
= \frac {\lambda}{z^\beta} w, \hbox{  } \forall (y, z)\in {\mathbf R}^{n-1}_y\times (0, \infty).
\end{equation}

Let us introduce also the time-dependent function:
\begin{equation}\label{W}W(t, y, z):=e^{{\bf i}\frac {\lambda t}{z^\beta}} w(y,z),
 \hbox{ } \forall (t, y, z)\in {\mathbf R}_t\times {\mathbf R}^{n-1}_y 
 \times (0, \infty),
\end{equation}
and let us compute now the partial differential equation satisfied by 
$W(t, y, z)$.

\vspace{0.2cm}

Elementary computations imply that
\begin{align*}
\partial_z^2 W =  e^{{\bf i}\frac{\lambda t}{z^\beta}}\left[
\frac{\beta(\beta+1){\bf i} \lambda t}{z^{\beta+2}} 
v\left (\frac{y}{\sqrt{z^\beta}} \right)
-\frac{\beta^2\lambda^2 t^2}{z^{2\beta+2}} 
v\left (\frac{y}{\sqrt{z^\beta}}\right )+
\frac{\beta^2{\bf i}\lambda t}{2z^{\frac32 \beta+2}} 
y. \nabla_y v \left(\frac {y}{\sqrt{z^\beta}}\right) \right. \\
\left. +\frac {\beta(\beta+2)}{4z^{\frac{\beta}2+2}} 
y. \nabla_y v \left(\frac {y}{\sqrt{z^\beta}}\right)+
 \frac {\beta^2{\bf i\lambda  }t}{2z^{\frac32 \beta+2}} 
y. \nabla_y v \left(\frac { y}{\sqrt{z^\beta}}\right)
+\frac {\beta^2} {4z^{\beta+2}} 
y. D^2_y v \left(\frac {y}{\sqrt{z^\beta}}\right).y \right]
\end{align*}
where $D^2_yv=\left (\frac{\partial^2 v }
{\partial y_i \partial y_j} \right)_{i,j=1,...,n-1}$.

\vspace{0.2cm}

Notice that by introducing the functions
\begin{equation}\label{GH}G(y)=y.\nabla_y v(y) 
\hspace{0.2cm} \hbox{ and } \hspace{0.2cm}
H(y)=y.D^2_yv(y).y, \hspace{0.2cm}
\forall y\in {\mathbf R}^{n-1}_y,
\end{equation}
the previous identity can be written as follows:

\begin{align*}
e^{-{\bf i}\frac {\lambda t}{z^\beta}} \partial_z^2 W
= \ 
&v\left (\frac {y}{\sqrt{z^\beta}}\right )\left(
\frac{\beta(\beta+1){\bf i}\lambda t}{z^{\beta+2}} - \frac{\beta^2\lambda^2 t^2}{z^{2\beta+2}}
\right ) \\
&+ \  
G \left (\frac {y}{\sqrt{z^\beta}}\right )
\left ( \frac{\beta^2{\bf i} \lambda t}{z^{\beta+2}} +\frac{\beta(\beta+2)}{4z^2}\right)
+\frac {\beta^2}{4 z^2} H\left (\frac { y}{\sqrt{z^\beta}}\right ).
\end{align*}

\vspace{0.2cm}

By using \eqref{rescaling} and the definition of $W$
(see \eqref{W}), we get 

\begin{align*}
{\bf i}\partial_t W  - \Delta_x W &= e^{{\bf i}\frac {\lambda t}{z^\beta}}
\left [-\frac {\lambda}{z^\beta} w - 
\Delta_y w - e^{-{\bf i}\frac{\lambda t}{z^\beta}} \partial_{z}^2 W \right ]\\
&= - \frac{1}{z^{2\beta}} \left ( \sum_{i,j=1}^{n-1} 
a_{ij} y_iy_j \right)e^{{\bf i}\frac {\lambda t}{z^\beta}} w- \partial_{z}^2 W
\end{align*}
that can be written in the following way
\begin{equation}\label{schr1}
{\bf i}\partial_t W - \Delta_x W +  
\frac{1}{z^{2\beta}} \left ( \sum_{i,j=1}^{n-1} 
a_{ij} y_iy_j \right) W= F , \forall (t, y, z) 
\in {\mathbf R}_t\times {\mathbf R}^{n-1}_y\times (0, \infty)
\end{equation}
where
\begin{align} \label{F} 
e^{-{\bf i}\frac {\lambda t}{z^\beta}} F (t, y, z) 
=\ \frac{1}{z^2} &\left[ v\left (\frac {y}{\sqrt{z^\beta}}\right )\left(
-\frac{\beta(\beta+1){\bf i} \lambda t}{z^\beta} + \frac{\beta^2\lambda^2 t^2}{z^{2\beta}}\right)\right. \notag \\
&\left. -\ G \left (\frac {y}{\sqrt{z^\beta}}\right )\left ( 
\frac{\beta^2{\bf i}\lambda t}{z^\beta} +\frac{\beta(\beta+2)}{4}\right ) 
-\frac {\beta^2}{4} H\left (\frac {y}{\sqrt{z^\beta}}\right )\right]. 
\end{align}

\vspace{0.2cm}

Let us introduce now the real valued cut-off function
$\phi(z) \in C^\infty
({\mathbf R}, [0,1])$
with the following properties:

\vspace{0.1cm}



\begin{enumerate}
\item $\phi(z)=0, \hbox{ for } |z|>1$,
\item  $\phi(z)=1, \hbox{ for } |z|<\frac 12$. 
\end{enumerate}

\vspace{0.1cm}

For a fixed parameter $\gamma \in (0, 1)$
we can truncate the function $W$ 
introduced in \eqref{W} in the following way:
$$W_{R}(t, y, z)=
W (t,y,z) 
\phi\left (\frac{z-R}{R^{\gamma}} \right)
\phi \left (\frac{|y|^2}{z^2} \right ),
\hbox{ }(t, y, z)\in {\mathbf R}_t 
\times {\mathbf R}^{n-1}_y \times {\mathbf R}_z.$$

\vspace{0.1cm}

Since now on we shall use the following notations:
\begin{align}\label{notations}
\phi_R &=\phi\left (\frac{z-R}{R^\gamma} \right),&
\phi_R' &=\frac 1{R^\gamma} \phi'\left (\frac{z-R}{R^\gamma} \right),&
\phi_R'' &=\frac 1{R^{2\gamma}}\phi''\left (\frac {z-R}{R^\gamma} \right), \\
\phi &= \phi\left(\frac{|y|^2}{z^2} \right),& 
\phi' &=\phi'\left(\frac{|y|^2}{z^2} \right),& 
\phi'' &=\phi''\left(\frac{|y|^2}{z^2} \right). \notag
\end{align}

We can write now the partial differential 
equation satisfied by $W_{R}$:

\begin{align*}
{\bf i}\partial_t W_{R} - &\Delta_x W_{R} 
\ =\
{\bf i}\phi_R  
\phi\partial_t 
W - \phi_R \Delta_y
(W \phi) 
- \partial_z^2 (W \phi_R 
\phi)
\\
=\ &\phi_R \phi ({\bf i}\partial_t 
W  - \Delta_y W  - 
\partial_z^2 W ) 
-\frac{4}{ z^2} \phi_R  
\phi'  y. 
\nabla_y W - W \partial_z^2 \left (\phi_R 
\phi \right )
\\
&- W  \phi_R \left (\frac {2(n-1)}
{ z^2} \phi' + \frac {4 |y|^2}
{ z^4}\phi''\right )
- 2 \partial_z W \left (
\phi_R' \phi 
-\frac{2|y|^2}{z^3} \phi_R \phi' \right ).
\end{align*}

\vspace{0.1cm}

By combining this identity with \eqref{schr1} 
and introducing the function 
\begin{equation*}
\Gamma(y, z)= G\left ( \frac y{\sqrt {z^\beta} }\right), \forall 
(y, z)\in {\mathbf R}_y^{n-1}\times (0, \infty),
\end{equation*}
where $G(y)$ is the function defined in \eqref{GH},
we get:

\begin{align*}
{\bf i}\partial_t W_{ R} - \Delta_x W_{R}
= &- \frac{1}{z^{2\beta}}\left (\sum_{i,j=1}^{n-1} 
a_{ij} y_i y_j \right) \phi_R \phi W + 
\phi_R \phi F
-\frac{4 e^{\frac{{\bf i}\lambda t}{z^\beta}}}{z^2} \phi_R 
\phi' \Gamma \\
&- W  \phi_R 
\left (\frac {2(n-1)}{z^2} \phi' + \frac {4 |y|^2}{ z^4}\phi'' \right) \\
&- W  \left (\phi_R''
\phi 
- \frac{4 |y|^2}{z^3} \phi_R' 
 \phi'
\right .   \left. +\frac{6  |y|^2}{ z^4} \phi_R 
\phi'
+ \frac{4|y|^4}{ z^6} \phi_R  
\phi''\right ) \\
&+2 \left (\frac {\beta{\bf i}\lambda t}{z^{\beta+1}}W 
+ \frac{\beta e^{\frac{{\bf i}\lambda t}
{z^\beta}}}{2z}\Gamma  \right)  \left (\phi_R'
\phi 
-\frac{2 |y|^2}{z^3} \phi_R 
\phi'\right ).
\end{align*}

\vspace{0.2cm}

Thus finally we have

\vspace{0.2cm}

\begin{equation}\label{CP1}
\begin{cases}
{\bf i}\partial_t W_{R} - \Delta_x 
W_{R}+ \frac{1}{z^{2\beta}}\left (\sum_{i,j=1}^{n-1}a_{ij}y_i y_j \right) 
W_{R} = F_{R}, \\ 
W_{R}(0, y, z)= f_{R},
\end{cases}
\end{equation}
where
\begin{equation}\label{fR}
f_{R}(y,z)=\phi_R 
\phi w,\ \forall (y, z)\in {\mathbf R}_y^{n-1}\times {\mathbf R}_z,
\end{equation}

\begin{equation}\label{FR}F_{R}(t, y, z)=\phi_R \phi
F + G_{R},\ \forall 
(t, y, z)\in {\mathbf R}_t
\times{\mathbf R}_y^{n-1}\times {\mathbf R}_z\end{equation}

and

\begin{align}
\label{GRLe} e^{-\frac{{\bf i} \lambda t}{z^\beta}}G_{R}(t, y,z)
= -w\bigg(&\frac {2(n-1) }{ z^2}\phi_R  \phi' 
+ 
\frac {4 |y|^2}{ z^4}\phi_R \phi''
+
\phi_R''\phi 
-
\frac{4 |y|^2}{z^3} \phi_R'  \phi'  \notag \\
&\left. + \frac{6 |y|^2}{ z^4} \phi_R \phi'
+ 
\frac{4|y|^4}{ z^6} \phi_R  \phi''
-
\frac{2\beta{\bf i} \lambda t}{z^{\beta+1}}\phi_R'\phi 
+
\frac{4\beta{\bf i} \lambda t  |y|^2}{z^{\beta+4}} \phi_R \phi'
\right) \\
 - \Gamma \bigg (
&\frac{4 }{z^2} \phi_R \phi' 
-
\frac{\beta}{z} \phi_R' \phi
+
\frac{2\beta|y|^2}{ z^4} \phi_R \phi'
\bigg ). \notag
\end{align}

The equation \eqref{CP1} and subsequent formulas for $f_R$, $F_R$, and
$G_R$ hold for every 
$(t, y, z)\in {\mathbf R}_t
\times{\mathbf R}_y^{n-1}\times {\mathbf R}_z$. 
The functions $f_R, W_R, F_R$ introduced
in this section will 
be very important in next section where we prove theorem \ref{main}.

\section{Proof of Theorem \ref{main}}
\label{sectionthm}

In  next lemma we shall assume
that the functions $f_R, W_R, F_R, $
are the ones constructed in the previous section starting from a 
fixed couple $\lambda, v(y)$ that satisfies \eqref{eigen}.
Notice that  $f_R, W_R, F_R,$ depend also on the parameter 
$\gamma$.

We shall also assume that the matrix 
$\left (a_{ij} \right)_{i,j=1,...,n-1}$
that appears in \eqref{eigen} is a fixed one.
We finally recall that 
we are using notations 
\eqref{notations} and that $\beta=1+ \frac \sigma 2$, where $\sigma$ is the same constant 
that appears in the assumptions
of theorem \ref{main}.
\vspace{0.1cm}

\begin{lem}
The following estimates are satisfied 
for any $n\geq 2$, $1\leq p, q<\infty$, and $0 < \gamma < 1$:
\begin{equation}\label{cauchy}\|f_{R}\|_{L^2_x}
\leq C R^\frac{(n-1)\beta+2\gamma}{4},  \hspace{0.2cm} \forall R>2,
\end{equation}
\begin{equation}\label{stri}\|W_{R}\|_{L^p_T L^q_x}
\geq C T^\frac 1p R^\frac{(n-1)\beta+2\gamma}{2q}, \hspace{0.2cm} 
\forall R>2, T>0,
\end{equation}
\begin{equation}\label{FconR}
\left \|F_R \right \|_{L^p_TL^q_x}
\leq
C  T^\frac 1p R^{\frac{((n-1)\beta+2\gamma)}{2q}}Max
\{ R^{-2\gamma}, T^{2}R^{-(2\beta+2)}\}, \hspace{0.2cm} \forall R>2, T>0,
\end{equation}
where $C=C(q, \gamma)>0$.

\noindent In particular if $$\frac 2p + \frac nq=\frac n2,$$
and $\alpha>0$ is any fixed number,
then
\begin{equation}\label{Wf}
\frac{\|W_{R}\|_{L^p((0,R^\alpha);L^q_x)}}
{\|f_{R}\|_{L^2_x}}
\geq C R^{\frac{\alpha n-((n-1)\beta+2\gamma)}{np}},\hspace{0.2cm} \forall R>2,
\end{equation}
and
\begin{equation}\label{WF}
\frac{\|W_{ R}\|_{L^{p}((0,R^\alpha);L^{q}_x)}}
{\|F_{R}\|_{L^{p'}((0, R^\alpha);L^{q'}_x)}}
\geq C R^{\kappa}, \hspace{0.2cm} \forall R>2,
\end{equation}
where
\begin{equation}\label{kappa}
\kappa= \kappa(n, \gamma, \alpha, p) \end{equation}
$$
=2\left(\frac{n\alpha -((n-1)\beta+2\gamma)}{np}\right) 
+ Min \left \{ 2\gamma - \alpha, 2\beta + 2 - 3\alpha \right \},
$$
and the constants $C>0$ do not depend on $R >2$.
\end{lem}

\noindent{\bf Proof.}
Assume that 
$L(y)$
is a given function,
then we introduce the new functions
\begin{equation}\label{lambda}\Lambda(y, z)= 
L\left (\frac{y}{\sqrt{z^\beta}}\right),
\Omega(y, z)=\left(\frac{|y|}{\sqrt{z^\beta}}\right)^2 
L\left (\frac{y}{\sqrt{z^\beta}}\right)
\end{equation}
\begin{equation*}
\hbox{ and } \Psi(y, z)= 
\left( \frac{|y|}{\sqrt{z^\beta}}\right)^4 L\left (\frac{y}
{\sqrt{z^\beta}}\right), \forall (y, z)\in {\mathbf R}_y^{n-1}\times
(0,\infty).
\end{equation*}

We claim that if
$L$ is a nontrivial function such that 
$L\in \cap_{p=1}^\infty L^p (\mathbf R^{n-1}_y)$, 
then the following estimate holds to be true:
\begin{equation}\label{thetaphi}
c R^\frac{(n-1)\beta+2\gamma}{2q}
\leq \left \|\Lambda
\phi_R \phi \right \|_{L^q_x} \leq
C  R^\frac{(n-1)\beta+2\gamma}{2q}, \hspace{0.2cm} \forall R>1,
\end{equation}
where $c=c(q,L)>0$ and $C=C(q,L)>0$, while $\phi_R$ and $\phi$ are the 
functions introduced in section \ref{preliminary}.

\vspace{0.1cm}

Notice that due to the properties of the functions
$\phi_R, \phi$
we can write the following chain of inequalities:
\begin{multline*}
\int_{R-\frac{R^\gamma}2}^{R+\frac{R^\gamma}2}dz 
\int_{|y|<\frac {\sqrt {z^{2-\beta}}}{\sqrt 2}}
 |L ( y )|^q z^\frac{(n-1)\beta}2 dy
= \int_{R-\frac{R^\gamma}2}^{R+\frac{R^\gamma}2} dz 
\int_{|y|<\frac{z}{\sqrt 2}}
\left |\Lambda \right |^q dy \\
\leq \int_{{\mathbf R}^n} 
\left |\Lambda\phi_R \phi \right |^q dy dz \\
\leq
\int_{R-{R^\gamma}}^{R+{R^\gamma}} dz \int_{|y|<z} 
\left |\Lambda\right |^qdy
= \int_{R-{R^\gamma}}^{R+{R^\gamma}} dz \int_{|y|< \sqrt {z^{2-\beta}}}
 |L ( y )|^q z^\frac{(n-1)\beta}2 dy,
\end{multline*}
that implies easily \eqref{thetaphi}.

\vspace{0.1cm}

Notice that with the same argument we can prove:
\begin{equation}\label{thetaphiprimo}
c R^\frac{(n-1)\beta+2\gamma}{2q}
\leq \left \|\Lambda
\phi_R \phi' \right \|_{L^q_x} \leq
C R^\frac{(n-1)\beta+2\gamma}{2q},
\end{equation}

\begin{equation}\label{thetaprimophi}
c  R^{\frac{(n-1)\beta+2\gamma}{2q}-\gamma}
\leq \left \|\Lambda
\phi_R' \phi\right \|_{L^q_x} \leq
C R^{\frac{(n-1)\beta+2\gamma}{2q}-\gamma},
\end{equation}
\begin{equation}\label{thetaprimophiprimo}
c R^{\frac{(n-1)\beta+2\gamma}{2q}-\gamma}
\leq \left \|  \Lambda
\phi_R'\phi' \right \|_{L^q_x} \leq
C R^{\frac{(n-1)\beta+2\gamma}{2q}-\gamma},
\end{equation}
\begin{equation}\label{thetaphisecondo}
c  R^\frac{(n-1)\beta+2\gamma}{2q}
\leq \left \|\Lambda
\phi_R \phi''\right \|_{L^q_x} \leq
C R^\frac{(n-1)\beta+2\gamma}{2q},
\end{equation}
\begin{equation}\label{thetasecondophi}
c  R^{\frac{(n-1)\beta+2\gamma}{2q}-2\gamma}
\leq  \left \|\Lambda
\phi_R''\phi \right \|_{L^q_x} \leq
C  R^{\frac{(n-1)\beta+2\gamma}{2q}-2\gamma}.
\end{equation}

\vspace{0.2cm}

{\em Proof of \eqref{cauchy}
and \eqref{stri}}

\vspace{0.2cm}

They follow easily from  
\eqref{thetaphi} 
where we choose $L(y)=v(y)$.

\vspace{0.2cm}

{\em Proof of \eqref{FconR}}

\vspace{0.2cm}

Due to \eqref{FR} it is easy to see that 
\eqref{FconR} follows from the following estimates:
\begin{equation}\label{Fepsilondelta23}
\left \|\phi_R \phi F \right \|_{L^p_TL^q_x}
\leq
C  T^\frac 1p R^{\frac{((n-1)\beta+2\gamma)}{2q}}Max
\{  R^{-2}, T^{2}R^{-(2\beta+2)}, \},
\end{equation}
and
\begin{equation}\label{Gepsilon}
\|G_{R}\|_{L^p_TL^q_x}
\leq C T^\frac 1p 
R^{\frac {((n-1)\beta+2\gamma)}{2q}} Max\{ 
R^{-2\gamma}, R^{-(6-2\beta)}, TR^{-(\beta+1+\gamma)}, T R^{-4}\}
\end{equation}
that we are going to prove, under the imposed conditions that $0<\gamma < 1$
and $1\le \beta < 2$ (recall that this range of $\beta$ corresponds to homogeneity
of the potential of order $0 \le \sigma < 2$.)  Some of the possible maximizers
can also be removed from consideration by observing geometric progressions.
To give one example, $TR^{-(\beta+1+\gamma)}$ is the geometric mean of
$R^{-2\gamma}$ and $T^2R^{-(2\beta+2)}$; therefore it must be intermediate in
size relative to the other two quantities.

\vspace{0.1cm}

Looking at the structure of $F$ (see \eqref{F}),
it is easy to see that 
in order to deduce \eqref{Fepsilondelta23} it is
sufficient to estimate
the norms of functions of the following type:
$$\frac{ t}{z^{\beta+2}}\Lambda \phi_R
\phi , 
\frac{t^2}{z^{2\beta+2}}
\Lambda \phi_R\phi, 
\frac{1}{z^2} \Lambda\phi_R \phi$$
where $\Lambda$ is defined as in
\eqref{lambda} and $L$ may change in different expressions.

\vspace{0.2cm}

Let us consider the first term.
Notice that due to the localization 
properties of the function $\phi_R(z)$,
we have that $\frac{t}{z^{\beta+2}}\leq \frac{C t}{R^{\beta+2}}$
on the support of the function $ \frac{t}{z^{\beta+2}}\Lambda \phi_R
\phi $,
thus with simple computations we get
\begin{equation} \label{Festimate1}
\left \|
\frac{t}{z^{\beta+2}}\Lambda \phi_R 
\phi  \right \|_{L^p_TL^q_x}\end{equation}$$
\leq C T^{1+\frac 1p}\frac{1}{R^{\beta+2}}
\left \|\Lambda \phi_R 
\phi \right \|_{L^q_x}
\leq C T^{1+\frac 1p} R^{\frac{(n-1)\beta+2\gamma}{2q}-(\beta+2)},
$$
where we have used \eqref{thetaphi}
in the last estimate.

With similar arguments we can deduce that
$$\left \|
\frac{ t^2}{z^{2\beta+2}}\Lambda \phi_R 
\phi \right \|_{L^p_TL^q_x}
\leq 
C T^{2+\frac 1p} R^{\frac{(n-1)\beta+2\gamma}{2q}-(2\beta+2)}
\qquad \hbox  { and }$$$$
\left \|
\frac{1}{z^2}\Lambda \phi_R
\phi \right \|_{L^p_TL^q_x}
\leq 
C T^{\frac 1p} R^{\frac{(n-1)\beta+2\gamma}{2q}-2}.$$
Note again that the norm estimate in \eqref{Festimate1} is the geometric mean 
of the two estimates above, so it cannot be the largest of the three
derived quantities.  
With this observation \eqref{Fepsilondelta23} is proved.

\vspace{0.1cm}

Looking at the structure of $G_{ R}$
(see \eqref{GRLe}) it is easy to see that \eqref{Gepsilon} 
comes from the following estimates
(whose proof follows as above by combining 
the localization properties of $\phi_R(z)$
with \eqref{thetaphiprimo}, \eqref{thetaprimophi},
\eqref{thetaprimophiprimo}, \eqref{thetaphisecondo}, 
\eqref{thetasecondophi})
where $\Lambda, \Omega$ and $\Psi$  
are the functions defined in \eqref{lambda}
and may depend on different $L$ in different expressions:

$$\left \|
\frac{1}{ z^2}\Lambda \phi_R 
\phi' \right \|_{L^p_TL^q_x}
\leq 
C T^{\frac 1p}  R^{\frac{(n-1)\beta+2\gamma}{2q}-2};$$

$$\left \|
\frac{|y|^2}{ z^4}\Lambda \phi_R
\phi'' \right \|_{L^p_TL^q_x}=
\left \|
\frac{1}{ z^{4-\beta}}\Omega \phi_R 
\phi'' \right \|_{L^p_TL^q_x}
\leq 
C T^{\frac 1p} R^{\frac{(n-1)\beta+2\gamma}{2q}-(4-\beta)};$$

$$\left \|
\Lambda \phi_R''
\phi  \right \|_{L^p_TL^q_x}
\leq 
C T^{\frac 1p} R^{\frac{(n-1)\beta+2\gamma}{2q}-2\gamma};$$

$$\left \|
\frac{|y|^2}{z^3}\Lambda \phi_R'
\phi' \right \|_{L^p_TL^q_x}=
\left \|
\frac{1}{ z^{3-\beta}}\Omega \phi_R'
\phi' \right \|_{L^p_TL^q_x}
\leq 
C T^{\frac 1p}  R^{\frac{(n-1)\beta+2\gamma}{2q}-(3-\beta+\gamma)};$$

$$\left \|
\frac{|y|^2}{ z^4}\Lambda \phi_R 
\phi' \right \|_{L^p_TL^q_x}=
\left \|
\frac{1}{z^{4-\beta}}\Omega \phi_R
\phi' \right \|_{L^p_TL^q_x}
\leq 
C T^{\frac 1p} R^{\frac{(n-1)\beta+2\gamma}{2q}-(4-\beta)};$$

$$\left \|
\frac{|y|^4}{ z^6}\Lambda \phi_R 
\phi''  \right \|_{L^p_TL^q_x}=
\left \|
\frac{1}{z^{6-2\beta}}\Psi \phi_R
\phi'' \right \|_{L^p_TL^q_x}
\leq 
C T^{\frac 1p}  R^{\frac{(n-1)\beta+2\gamma}{2q}-(6-2\beta)};$$

$$\left \|
\frac{ t}{z^{\beta+1}}\Lambda \phi_R'
\phi \right \|_{L^p_TL^q_x}
\leq 
C T^{1+\frac 1p} R^{\frac{(n-1)\beta+2\gamma}{2q}-(\beta+1+\gamma)};$$

$$\left \|
\frac{ t |y|^2}{
z^{\beta+4}}\Lambda \phi_R
\phi' \right \|_{L^p_TL^q_x}=
\left \|
\frac{t}{z^4}\Omega \phi_R 
\phi' \right \|_{L^p_TL^q_x}
\leq 
C T^{1+\frac 1p} R^{\frac{(n-1)\beta+2\gamma}{2q}-4};$$

$$\left \|
\frac{1}{z}\Lambda \phi_R'
\phi \right \|_{L^p_TL^q_x}
\leq 
C T^{\frac 1p} R^{\frac{(n-1)\beta+2\gamma}{2q}-(1+\gamma)}.$$

The estimate \eqref{Gepsilon} now follows easily.

\vspace{0.2cm}

{\em Proof of \eqref{Wf} and \eqref{WF}}

\vspace{0.2cm}

They follow with elementary computations from
\eqref{cauchy}, \eqref{stri}
and \eqref{FconR} where we choose $T=R^\alpha$. 

\vspace{0.2cm}

\hfill$\Box$

\vspace{0.3cm}

\noindent{\bf Proof of theorem \ref{main}.}
We can assume without loss of generality that the minimum 
of the restriction of $V(x)$ on the sphere ${\mathbf S}^{n-1}$
is achieved at 
the point $(0,...0,1)\in {\mathbf R}_y^{n-1} \times {\mathbf R}_z$,
and that $V(x) = 0$ at this point.
Therefore
$$V(y, 1)= \sum_{i,j=1}^{n-1} a_{ij}y_i y_j
+ R(y), \hspace{0.2cm} \forall y\in {\mathbf R}_y^{n-1},$$
where the matrix $\left (a_{ij}\right)_{i,j=1,...,n-1}$
is the Hessian of the function $V(y, 1)$ at the point $y=0$,
which
is positive definite since we are assuming that
$V$ is a  generalized Morse type function and that $V(0,1) = 0$.

Moreover we have
$$\limsup_{|y|\rightarrow 0} |R(y)| |y|^{-3}<\infty,$$
and in particular the following pointwise estimate holds:
\begin{equation}\label{pointwise}|R(y)|\leq C |y|^3, 
\hspace{0.2cm} \forall y\in {\mathbf R}^{n-1}_y
\hbox{ s.t. } |y|<1.
\end{equation}

Notice that since $V(x)$ is homogeneous of order $-\sigma$ we have that
$$V(y, z)= z^{-\sigma} V\left (\frac y z, 1\right)=
\frac 1{z^{2\beta}}\sum_{i,j=1}^{n-1} a_{ij}y_i y_j
+ z^{-\sigma} R\left (\frac yz \right), \hspace{0.2cm} \forall 
(y, z)\in {\mathbf R}_y^{n-1}\times (0, \infty),$$
(we have used $\beta=1+\frac \sigma2$)
and due to \eqref{pointwise}
we have that 
\begin{equation}\label{point}
\left |R\left (\frac yz \right)\right |\leq C \frac{|y|^3}{z^3},
\hbox{ provided that } |y|< |z|.
\end{equation}

We have now a well-defined positive definite symmetric 
matrix $\left (a_{ij}\right)_{i,j=1,...,n-1}$ and
we can select a couple
$(\lambda, v(y))$ that satisfies the eigenvalue problem \eqref{eigen}.

\vspace{0.1cm} 
Following the previous section 
we can construct starting from 
these fixed $(\lambda, v(y))$ the family of functions
$f_R, F_R, W_R$ that we shall use below.

\vspace{0.1cm}

Notice that due the cut-off property 
of the function $\phi$ and to \eqref{point} 
we get easily the following estimate: 
\begin{equation*}
\left \|z^{-\sigma} R\left ( \frac yz 
\right) W_R \right \|_{L^{p'}((0, T); L^{q'}_x)} 
\leq C T^{\frac 1{p'}} \left \|\frac{|y|^3}{|z|^{2\beta+1}} w \phi_R 
\phi \right \|_{L^{q'}_x}
$$$$\leq CT^{\frac 1{p'}}R^{-(\frac{\beta}2+1)} \left \|
M \phi_R
\phi \right \|_{L^{q'}_x}
\leq CT^{\frac 1{p'}}R^{-(\frac{\beta}2+1)+\frac{(n-1)\beta+2\gamma}{2q'}}, 
\hspace{0.1cm} \forall R>1, 0<\gamma<1
\end{equation*}
where $\beta = 1 + \frac{\sigma}2$ as always, 
$M(y,z)=\left ( \frac{|y|}{\sqrt {z^\beta} }\right )^3 
v\left (\frac y{\sqrt{z^\beta}}\right )$
and we have used \eqref{thetaphi} with $L(y)=|y|^3 v(y)$ at the last step.

\vspace{0.1cm}

As a by product of this inequality we get
\begin{equation}\label{ratio}
\left \|z^{-\sigma} R\left (\frac yz \right)
W_R\right \|_{L^{p'}((0, R^\alpha);L^{q'}_x)}
\leq {R^{\frac{\alpha}{p'}-(\frac{\beta}2+1)+ \frac{(n-1)\beta+2\gamma}{2q'}}},
\end{equation}
for any $\alpha>0$.
In particular if we assume that $\frac 2p + \frac nq=\frac n2$
and we use \eqref{stri} with $T=R^\alpha$, 
then it is easy to deduce that
\begin{equation}\label{ben}
\frac{\|W_R\|_{L^p((0, R^\alpha);L^q_x)}}
{\|z^{-\sigma} R\left(\frac yz \right)W_R\|_{L^{p'}((0, R^\alpha);L^{q'}_x)} }
\geq C R^{2\left(\frac{n\alpha - ((n-1)\beta+2\gamma)}{np}\right)+ 
(\frac{\beta}2+1-\alpha)}.
\end{equation}

Let us consider
now the following auxiliary Cauchy problems 
with non trivial forcing term:
\begin{equation}
\label{aux}
\begin{cases}{\bf i} \partial_t u_{R} - \Delta_x u_{ R}
+ \left [ \frac{1}{z^{2\beta}}\left (
\sum_{i,j=1}^{n-1} a_{ij}y_iy_j \right) 
+ z^{-\sigma} R\left ( \frac yz \right)\right ]u_R = \tilde F_{ R},\\
u_{R}(0, y, z)=f_{R}
\end{cases}\end{equation}
where
$$\tilde F_{R}(t, y, z):=\chi_{(0, R^\alpha)}(t) \left [F_{R}
+ z^{-\sigma} R\left (\frac yz  \right)W_R\right]$$ 
 and $f_R,F_R$ are given by \eqref{fR} and \eqref{FR}.

\vspace{0.1cm}

Let us notice that due to \eqref{WF} and \eqref{ben} 
we get
\begin{equation}\label{from}
\frac{\|W_R\|_{L^p((0, R^\alpha);L^q_x)}}
{\|\tilde F_R\|_{L^{p'}((0, R^\alpha);L^{q'}_x)} }
\geq C R^\delta, \forall 0<\gamma<1,
\end{equation}
where 
\begin{equation}\label{beta}
\delta=\delta(n, \gamma, \alpha, p)=2 \left(\frac{n\alpha -((n-1)\beta+2\gamma)}{np}
\right) \end{equation}
$$+ Min \left \{ 2\gamma - \alpha, 2\beta + 2 - 3\alpha, 
\frac{\beta}2 +1 - \alpha\right \}.
$$

Also due to \eqref{CP1} we have:
$$u_{R}(t, x)= W_{R}(t, x), \forall (t, x)\in 
(0, R^\alpha)\times {\mathbf R}_x^n,$$
then
$$\|u_{ R}\|_{L^p_t L^q_x}
\geq \|u_{R}\|_{L^p((0, R^\alpha); L^q_x)}= 
\|W_{ R}\|_{L^p((0, R^\alpha); L^q_x)},$$
from which it follows that
\begin{equation}\label{inequa}\frac
{\|u_{R}\|_{L^p_t L^q_x}}{\|f_{ R}\|_{L^2_x}
+ \|\tilde F_{R}\|_{L^{p'}_tL^{q'}_x}}
\ \geq\  \frac{\|W_{ R}\|_{L^p((0, R^\alpha); L^q_x)}}{\|f_{ R}\|_{L^2_x}
+ \|\tilde F_{ R}\|_{L^{p'}((0, R^\alpha);L^{q'}_x)}}.
\end{equation}

Notice that
\eqref{Wf} implies that for any $n\geq 2$ 
and for any $1\leq p < \infty$, $0<\gamma<1$,
$$\frac{\|W_{ R}\|_{L^p((0, R^\alpha); L^q_x)}}{\|f_{ R}\|_{L^2_x}}
\rightarrow \infty \hbox{ for } R\rightarrow \infty,$$
provided that $\alpha >\frac{(n-1)\beta+2\gamma}{n}$.

On the other hand the function $\delta (n,\gamma,\alpha, p)$, defined
in \eqref{beta}, varies
continuously with $\alpha$, and has the particular value
\begin{equation*}
\delta(n,\gamma, {\textstyle \frac{(n-1)\beta+2\gamma}{n}}, p) \end{equation*}
$$= Min \left \{
\frac{(n-1)(2\gamma-\beta)}{n}, \frac{(2-\beta)n + 3\beta -6\gamma}{n}, 
\frac{(2-\beta)n + 2\beta - 4\gamma}{2n}
\right\}.$$
This is strictly positive provided  $n \ge 2$ and 
$\frac{\beta}2 <\gamma < \frac{\beta}2 + \frac{(2-\beta)n}6$.
It is therefore possible to choose a value of $\gamma\in (0, 1)$ in this range,
then select 
$\alpha > \frac{(n-1)\beta+2\gamma}n$ so that
$\delta(n,\gamma,\alpha,p) > 0$.
By \eqref{from}, with these choices we get
$$\frac{\|W_{ R}\|_{L^p((0, R^\alpha); L^q_x)}}
{\|\tilde F_{ R}\|_{L^{p'}((0, R^\alpha);L^{q'}_x})}\rightarrow \infty
\hbox{ for } R\rightarrow \infty.$$

We can now deduce by \eqref{inequa} that
$$\frac{\|u_{ R}\|_{L^p_t L^q_x}}{\|f_{ R}\|_{L^2_x}
+ \|\tilde F_{ R}\|_{L^{p'}_tL^{q'}_x}}\rightarrow \infty \hbox{ for }
R\rightarrow \infty.
$$

Notice that with this last inequality
we have shown that 
estimates of the following type: 
$$\|u\|_{L^p_tL^q_x}\leq C \left (\|f\|_{L^2_x}
+ \|F\|_{L^{p'}_tL^{q'}_x} \right),$$
cannot be satisfied by the 
solutions of the following inhomogeneous Schr\"odinger equation
$$
\begin{cases}
{\bf i}\partial_t u - \Delta_x u+ V(x) u=F, (t, x)\in 
{\mathbf R}_t\times {\mathbf R}^n_x, \\
 u(0,x)=f(x),
\end{cases}
$$
when $V$ satisfies the assumptions of theorem \ref{main}. 
In order to disprove Strichartz estimates for the corresponding
homogeneous Schr\"odinger equation (i.e. 
for the Schr\"odinger equation with trivial forcing term)
in the case $p>2$, it is sufficient to combine the standard $TT^*$ argument
with a well-known result due to M. Christ and A. 
Kiselev (see \cite{CK}). 
As a consequence of this fact we can deduce that Strichartz estimates 
must be false also for $p=2$, otherwise by combining this estimate
with the trivial 
$L^\infty_tL^2_x$ estimate, we could get with an elementary 
interpolation argument the estimates also for $p>2$.

\hfill$\Box$


\begin{thebibliography}{CTCR81}

\bibitem{BRV} {\em J.A. Barcelo, A. Ruiz and L. Vega}
Some dispersive estimates
for Schr\"odinger equations with repulsive potentials.
Preprint. 

\bibitem{BPST1} {\em N. Burq, F. Planchon, J. Stalker and S. Tahvildar-Zadeh},
Strichartz estimates for the wave and Schr\"odinger
equations with the inverse-square potential,
J. Funct. Anal.,
vol. 203,
2003,
(2),
pp. 519-549.

\bibitem{BPST2} {\em N. Burq, F. Planchon, J. Stalker and S. Tahvildar-Zadeh},
Strichartz estimates for the wave and {S}chr\"odinger
equations with potentials of critical decay,
Indiana Univ. Math. J.,
vol. 53,
2004,
(6),
pp.1665-1680.

\bibitem{CK} {\em M. Christ and A. Kiselev},
Maximal functions associated to filtrations,
J. Funct. Anal.,
vol. 179,
2001,
(2),
pp. 409-425.

\bibitem{CS}{\em P. Constantin and J.C. Saut},
Local smoothing properties of dispersive equations,
J. Amer. Math. Soc.,
vol.(1),
1988,
(2),
pp. 413-439.

\bibitem{CoSa} {\em P. Constantin and J.C. Saut} 
Local smoothing properties of Schr\"odinger equations,
Indiana Univ. Math. J.,
vol. 38,
1989,
(3), pp. 791-810.




\bibitem{DPV} {\em P. D'Ancona, V. Pierfelice and N. Visciglia}
Some remarks on the
Schr\"odinger equation with potential in 
$L^r_t L^s_x$,
Math. Ann., vol. 333, 2005, pp.
271-290.


\bibitem{G}{\em M. Goldberg} Dispersive bounds for the 
three-dimensional Schr\"odinger equation with almost critical potentials,
to appear in Geom. and Funct. Anal.. 


\bibitem{H}{\em I. Herbst}, 
Spectral and scattering theory for {S}chr\"odinger operators
with potentials independent of $|x|$,
Amer. J. Math.,
vol. 113,
1991,
(3),
pp. 509-565.

\bibitem{HS}{\em I. Herbst and E. Skibsted}
Quantum scattering for potentials independent of $|x|$: 
asymptotic completeness for high and low energies,
Comm. Partial Differential Equations,
vol. 29, 2004, (3-4), pp. 547-610.
 

\bibitem{JSS}{\em  J.L. Journ\'ee, A. Soffer and C. Sogge},
Decay estimates for {S}chr\"odinger operators,
Comm. Pure Appl. Math.,
vol. 44,
1991,
(5),
 pp. 573-604.











\bibitem{KT} {\em M. Keel and T. Tao}, 
Endpoint Strichartz estimates,
Amer. J. Math.,
vol. 120,
1998,
(5),
pp. 955-980.

\bibitem{KPV}
{\em C. Kenig, G. Ponce and L. Vega} 
Small solutions for nonlinear Schr\"odinger equations,
Ann. Inst. Henri Poincar\'e Anal. Nonlin\'eaire,
vol. 10,
1993,
(3),
pp. 255-288.


\bibitem{PV1}{\em B. Perthame and L. Vega},
Energy concentration and {S}ommerfeld condition for
{H}elmholtz and {L}iouville equations,
C. R. Math. Acad. Sci. Paris,
vol. 337,
2003,
(9),
pp. 587-592.


\bibitem{PV2}{\em B. Perthame and L. Vega}
Energy decay and Sommerfield condition for
Helmholtz equation with variable index at infinity, 
Preprint.



\bibitem{RS} {\em I. Rodnianski and W. Schlag},
Time decay for solutions of Schr\"odinger equations with
rough and time-dependent potentials,
Invent. Math.,
vol. 155,
2004,
(3),
pp. 451-513.

\bibitem{RV1} {\em A. Ruiz and L. Vega},
On local regularity of Schr\"odinger equations,
 Internat. Math. Res. Notices,
1993,
(1),
pp. 13-27.



\bibitem{RV2} {\em A. Ruiz and L. Vega},
Local regularity of solutions to wave equations with
time-dependent potentials,
Duke Math. J.,
vol. 76,
1994,
(3),
pp. 913-940.


\bibitem{S}{\em W. Schlag}, Dispersive estimates for Schr\"odinger 
operators: a survey, Preprint.

\bibitem{Sj} {\em P. Sj{\"o}lin},
Regularity of solutions to the Schr\"odinger equation,
Duke Math. J., vol. 55, 1987, (3),
pp. 699-715.
   
\bibitem{ST}{\em G. Staffilani and D. Tataru,}
Strichartz estimates for a {S}chr\"odinger operator with
nonsmooth coefficients,
Comm. Partial Differential Equations,
vol. 27, 2002,
(7-8),
pp. 1337--1372.


\bibitem{V} {\em L. Vega},
Schr\"odinger equations: pointwise convergence to the initial
data,
Proc. Amer. Math. Soc.,
vol. 102,
1988,
(4),
pp. 874-878.

\bibitem{VV} {\em L. Vega and N. Visciglia},
On the local smoothing for the Schr\"odinger equation,
Preprint.


\end{thebibliography}
\end{document}